\documentclass[12pt, reqno]{amsart}
\usepackage{amsmath, amsthm, amscd, amsfonts, amssymb, graphicx, color}
\usepackage[bookmarksnumbered, colorlinks, plainpages]{hyperref}

\newtheorem{theorem}{Theorem}[section]

\newtheorem{proposition}[theorem]{Proposition}
\newtheorem{corollary}[theorem]{Corollary}
\theoremstyle{definition}
\newtheorem{definition}[theorem]{Definition}
\newtheorem{example}[theorem]{Example}

\theoremstyle{remark}
\newtheorem{remark}[theorem]{Remark}
\numberwithin{equation}{section}

\allowdisplaybreaks
\begin{document}

\title[$E$-$g$-frames]
{$E$-$g$-frames}

\author[H. Hedayatirad]{Hassan Hedayatirad}
\address{Hassan Hedayatirad \\ Department of Mathematics and Computer
Sciences, Hakim Sabzevari University, Sabzevar, P.O. Box 397, IRAN}
\email{hasan.hedayatirad@hsu.ac.ir; Hassan.hedayatirad67@gmail.com\rm }
\author[T.L. Shateri]{Tayebe Lal Shateri }
\address{Tayebe Lal Shateri \\ Department of Mathematics and Computer
Sciences, Hakim Sabzevari University, Sabzevar, P.O. Box 397, IRAN}
\email{ \rm  t.shateri@hsu.ac.ir; t.shateri@gmail.com}
\thanks{*The corresponding author:
t.shateri@hsu.ac.ir; t.shateri@gmail.com (Tayebe Lal Shateri)}
 \subjclass[2010] {Primary 42C15;
Secondary 54D55.} \keywords{$E$-$g$-frame, $E$-$g$-Bessel sequence, frame operator, direct sum of Hilbert spaces. }
 \maketitle

\begin{abstract}
In the present paper, we introduce the notion of $E$-$g$-frames for a separable Hilbert space $\mathcal H$, where $E$ is an invertible infinite matrix mapping on the Hilbert space $\mathop\oplus\limits_{n=1}^{\infty}\mathcal H_n$. We study some properties of $E$-$g$-frames. Also, we give a result concerning the perturbation of $E$-$g$-frames and then use it to construct $E$-$g$-frames in separable Hilbert spaces.
\vskip 3mm
\end{abstract}

\section{Introduction and Preliminaries}
The concept of frames has been introduced by Duffin and Schaeffer \cite{DS} to study some problems in non-harmonic Fourier series. Then Daubecheies, Grassman and Mayer \cite{DG} reintroduced and developed them. we refer readers to \cite{CH} for an introduction to frame theory in Hilbert spaces and its applications. Frames have very important and interesting properties which make them very useful in the characterization of function spaces, signal processing and many other fields such as image processing, data compressing, sampling theory and so on.

Various extensions of the frame theory have been investigated, several of them were contained in the theory of $g$-frames and $K$-$g$-frames (\cite{SUN,XIA,YAN}). Sun \cite{SUN} introduced $g$-frames and  $g$-Riesz bases as other generalized frames. He showed that oblique frames, pseudo-frames, and fusion frames are special cases of $g$-frames, also
from a $g$-frame, we may construct a frame for a complex Hilbert space. Xiao and et.al  \cite{XIA} introduced the concept of $K$-$g$-frames which extends the concepts of $K$-frames and $g$-frames. This fact caused several authors to study various aspects of $K$–$g$-frames. One of them is to get the methods of construction of $K$–$g$-frames. For instance, Du and Zhu \cite{ZHU} used the relation between positive operators and frame operators a $K$–$g$-frame to generate a new $K$–$g$-frame.

In this paper, by using some ideas from \cite{TAL} and \cite{JAV}, we present the definition of an $E$-$g$- frame in a separable Hilbert space and give some results of frames in the view of $E$-$g$-frames. $E$-$g$-frames in separable Hilbert spaces have some properties similar to those of frames, but not all the properties are similar. "Which properties of the frame may be extended to the $E$-$g$-frame?". We shall discuss this natural question.

Let $\mathcal X$ and $\mathcal Y$ be two sequence spaces and $E=(E_{n,k})_{n,k\geq 1}$ be an infinite matrix of real or complex numbers. We say that $E$ defines a matrix mapping from $\mathcal X$ into $\mathcal Y$, if for every sequence $x=\{x_n\}_{n=1}^{\infty}$ in $\mathcal X$, the sequnce $Ex=\{(Ex)_n\}_{n=1}^{\infty}$ is in $\mathcal Y$, where
$$(Ex)_n=\sum_{k=1}^{\infty}E_{n,k}x_k,\quad n=1,2,\cdots.$$

Throughout this paper, the space $\mathcal H$ denotes an infinite dimensional separable
Hilbert space and ${H_n}n\in\mathbb N$ is a sequence of closed subspaces of $\mathcal H$.  Let $\mathcal B(\mathcal H,\mathcal H_n)$ be the space of all bounded linear operators
from $\mathcal H$ into $\mathcal H_n$.\\
For each sequence $\{H_n\}_{n\in\mathbb N}$, we define the space $(\mathop\oplus\limits_{n=1}^{\infty}\mathcal H_n)_{l^2}$ by
$$(\mathop\oplus\limits_{n=1}^{\infty}\mathcal H_n)_{l^2}=\left\{\{f_n\}_{n=1}^{\infty}:f_n\in\mathcal H_n,\;\sum_{n=1}^{\infty}\|f_n\|^2<\infty \right\}.$$
With the inner product $\langle.,.\rangle$ defined by
$$\langle\{f_n\}_{n=1}^{\infty},\{g_n\}_{n=1}^{\infty}\rangle=\sum_{n=1}^{\infty}\langle f_n,g_n\rangle,$$
$(\mathop\oplus\limits_{n=1}^{\infty}\mathcal H_n)_{l^2}$ is a Hilbert space.

Recall that a sequence $\{\Lambda_n\in\mathcal B(\mathcal H,\mathcal H_n):n\in\mathbb N\}$ is said to be a $g$-frame for $\mathcal H$ if there exist positive real numbers $A$ and $B$ such that
\begin{equation}\label{def1}
A\left\|f\right\|^2\leq \sum_{n=1}^{\infty}\left\|\Lambda_nf\right\|^2\leq B\left\|f\right\|^2,
\end{equation}
for all $f\in \mathcal H$. The constants $A$ and $B$ are called lower and upper $g$-frame bounds. We recall the definition of an $E$-frame.
\begin{definition}\cite{TAL}
A sequence $\{f_k\}_{k=1}^{\infty}$ in a separable Hilbert spaces $\mathcal H$ is called an $E$-frame for $\mathcal H$ if there exist positive real numbers $A$ and $B$ such that
\begin{equation}\label{eq0}
A\left\|f\right\|^2\leq\left\|\left\{\left\langle f,\left(E\{f_j\}_{j=1}^{\infty}\right)_n\right\rangle\right\}_{n=1}^{\infty}\right\|_{\ell^2}^2\leq B\left\|f\right\|^2,
\end{equation}
for all $f\in\mathcal H$. the numbers $A$ and $B$ are called $E$-frame bounds. Inequality (\ref{eq0}) also can be written as
\begin{equation*}
A\left\|f\right\|^2\leq\sum_{n=1}^{\infty}\left|\left\langle f,\sum_{k=1}^{\infty}E_{n,k}f_k\right\rangle\right|^2\leq B\left\|f\right\|^2,
\end{equation*}
for all $f\in\mathcal H$.
\end{definition}
The pre $E$-frame operator $T$ is defined by
\begin{align*}
T:\ell^2\to\mathcal H,\quad T\{c_k\}_{k=1}^{\infty}=\sum_{k=1}^{\infty}c_k\left(E\{f_j\}_{j=1}^{\infty}\right)_k,
\end{align*}
is bounded. Its adjoint, the analysis operator, is given by
\begin{align*}
T^*:\mathcal H\to\ell^2,\quad T^*f=\left\{\left\langle f,\left(E\{f_j\}_{j=1}^{\infty}\right)_k\right\rangle\right\}_{k=1}^{\infty}.
\end{align*}
Composing $T$ and $T^*$, the $E$-frame operator
\begin{align*}
S=TT^*:\mathcal H\to\mathcal H,\quad Sf=\sum_{k=1}^{\infty}\left\langle f,\left(E\{f_j\}_{j=1}^{\infty}\right)_k\right\rangle\left(E\{f_j\}_{j=1}^{\infty}\right)_k.
\end{align*}
is obtained. The $E$-frame operator is bounded, invertible, self-adjoint and positive (see \cite{TAL}).
\section{$E$-$g$-frames}
In this section, we present the definition of an $E$-$g$-frame, then we give several construction methods of $E$-$g$-frames. More precisely, we give a result concerning perturbations of $E$-$g$-frames and then use it to construct $E$-$g$-frames in Hilbert spaces.
\begin{definition}
A sequence $\Lambda=\{\Lambda_n\}_{n\in\mathbb N}$ in $\mathcal B(\mathcal H,\mathcal H_n)$ is called an $E$-$g$-frame for $\mathcal H$ if there exist positive real numbers $A$ and $B$ such that
\begin{equation}\label{def1}
A\left\|f\right\|^2\leq \sum_{n=1}^{\infty}\left\|\sum_{k=1}^{\infty}E_{n,k}\Lambda_kf\right\|^2\leq B\left\|f\right\|^2,
\end{equation}
for all $f\in \mathcal H$. the numbers $A$ and $B$ are called $E$-$g$-frame bounds of $\{\Lambda_n\}_{n\in\mathbb N}$. An $E$-$g$-frame is called a tight frame if $A = B$, and is called a Parseval frame if $A = B = 1$. If only the right hand inequality of \ref{def1} holds, then we say that $\{\Lambda_n\}_{n\in\mathbb N}$ is an $E$-$g$-Bessel sequence.
\end{definition}
Note that if $\mathcal H_n=\mathcal H$ and $\Lambda_n=id_{\mathcal H}$, for every $n\in\mathbb N$, then an $E$-$g$-frame is an $E$-frame.

In the sequel, we give two examples of an $E$-$g$-frame.
\begin{example}\label{exam1}
Let ${\bf\Delta}=\left({\bf\Delta}_{n,k}\right)_{n,k\geq 1}$ denotes the matrix defined by
\begin{equation*}
{\bf\Delta}= \left\{
\begin{array}{rl}
    (-1)^{n-k}  &  n-1\leq k\leq n \\
    0  \qquad  & 1\leq k\leq n-1\; \text{or}\; k>n.
\end{array} \right.
\end{equation*}
The matrix form of $\bf\Delta$ is 
\begin{equation*}
\left(
\begin{matrix}
1  & 0 & 0 & 0 & 0 &\cdots \\
   -1  & 1 & 0 & 0 & 0 &\cdots \\
    0  & -1 & 1 & 0 & 0 &\cdots \\
    0  & 0 & -1 & 1 & 0 & \cdots \\
    0  & 0 & 0 & -1 & 1 &\cdots \\
    \vdots & \vdots & \vdots & \vdots & \ddots &  \ddots\\
   \end{matrix}
\right).
\end{equation*}
Then for a sequence $\{f_j\}_{j=1}^{\infty}$ in $\mathcal H$ we have
\begin{equation*}
\left({\bf\Delta}\{f_j\}_j\right)_n=\sum_{k=1}^{\infty}{\bf\Delta}_{n,k}f_k=f_n-f_{n-1}.
\end{equation*}
Let $\{f_n\}_{n=1}^{\infty}$ be an arbitrary frame in $\mathcal H$ with the frame bounds $A$ and $B$. Define $\Lambda_{f_n}$ 
be the functional induced by $f_n$ , i.e.,
$$\Lambda_{f_n}f=\left\langle f, f_n\right\rangle\qquad (f\in\mathcal H),$$
and we put $\Lambda_0f=0$. Then $\{\Lambda_n\}_n$ is a $\bf\Delta$-$g$-Bessel sequence for 
$\mathcal H$ with the $\bf\Delta$-Bessel bound  $4B$.
\begin{align*}
\sum_{n=1}^{\infty}\left|\sum_{k=1}^{\infty}{\bf\Delta}_{n,k}\Lambda_{f_k}f\right|^2&=\sum_{n=1}^{\infty}\left|\Lambda_{f_n} f-\Lambda_{f_{n-1}} f\right|^2\\
&=\sum_{n=1}^{\infty}\left|\left\langle f,f_n\right\rangle-\left\langle f,f_{n-1}\right\rangle\right|^2\\
&\leq\sum_{n=1}^{\infty}2\left(\left|\left\langle f,f_n\right\rangle\right|^2+\left|\left\langle f,f_{n-1}\right\rangle\right|^2\right)\\
&=2\left(\sum_{n=1}^{\infty}\left|\left\langle f,f_n\right\rangle\right|^2+\sum_{n=1}^{\infty}\left|\left\langle f,f_{n-1}\right\rangle\right|^2\right)\\
&=4\sum_{n=1}^{\infty}\left|\left\langle f,f_n\right\rangle\right|^2\\
&\leq 4B.
\end{align*}
\end{example}
\begin{example}\label{exam2}
With the same hypotheses in Example \ref{exam1}, the sequence \\$\{g_n\}_{n=1}^{\infty}=\{f_1,0,f_2,0,f_3,0,\cdots\}$ is a both a frame with the same frame bounds as $\{f_n\}_{n=1}^{\infty}$. Furthermore, $\{\Lambda_{g_n}\}_{n=1}^{\infty}$ is a $\bf\Delta$-$g$-frame for 
$\mathcal H$ with $\bf\Delta$-frame bounds  $A$ and $2B$. 
\end{example}
It is clear that the operator
\begin{equation*}
T_{\Lambda}:(\mathop\oplus\limits_{n=1}^{\infty}\mathcal H_n)_{l^2}\to\mathcal H,\quad T_{\Lambda}\{f_n\}_{n=1}^{\infty}
=\sum_{n=1}^{\infty}\sum_{k=1}^{\infty}\overline{E}_{n,k}\Lambda_k^*f_n,
\end{equation*}
is bounded. $T$ is called the pre $E$-$g$-frame operator. The analysis operator, is given by
\begin{equation*}
T_{\Lambda}^*:\mathcal H\to (\mathop\oplus\limits_{n=1}^{\infty}\mathcal H_n)_{l^2},\quad T_{\Lambda}^*f=\left\{\sum_{j=1}^{\infty}E_{n,j}\Lambda_jf\right\}_{n=1}^{\infty}.
\end{equation*}
The $E$-$g$-frame operator is defined as 
\begin{equation*}
S_{\Lambda}=T_{\Lambda}T_{\Lambda}^*:\mathcal H\to\mathcal H,\quad S_{\Lambda}f
=\sum_{n=1}^{\infty}\sum_{k=1}^{\infty}\sum_{j=1}^{\infty}\overline{E}_{n,k}E_{n,j}\Lambda_k^*\Lambda_jf.
\end{equation*}
We can easily prove the following result.
\begin{proposition}
Let $\Lambda=\{\Lambda_n\}_{n\in\mathbb N}$ be an $E$-$g$-frame with frame bounds $A$ and $B$. Then the $E$-$g$-frame operator $S_{\Lambda}$ is bounded, invertible, self-adjoint, and positive and $A\leq S_{\Lambda}\leq B$. 
\end{proposition}
\begin{remark}
For any $f\in\mathcal H$, we have
\begin{align*}
f=S_{\Lambda}S_{\Lambda}^{-1}=S_{\Lambda}^{-1}S_{\Lambda}f=\sum_{n=1}^{\infty}\sum_{k=1}^{\infty}\sum_{j=1}^{\infty}\overline{E}_{n,k}E_{n,j}\Lambda_k^*\Lambda_jS_{\Lambda}^{-1}f.
\end{align*}
The above $E$-$g$-frame decomposition shows that every element in $\mathcal H$ has a representation as an infinite linear combination of the $E$-$g$-frame elements.
\end{remark}
Let $\tilde{\Lambda}_n=\Lambda_nS_{\Lambda}^{-1}$, for each $n\in\mathbb N$, then  the above equalities become
\begin{align*}
f=\sum_{n=1}^{\infty}\sum_{k=1}^{\infty}\sum_{j=1}^{\infty}\overline{E}_{n,k}E_{n,j}\Lambda_k^*\tilde{\Lambda}_jf=\sum_{n=1}^{\infty}\sum_{k=1}^{\infty}\sum_{j=1}^{\infty}\overline{E}_{n,j}E_{n,k}\tilde{\Lambda}_j^*\Lambda_kf.
\end{align*}
We give the following result.
\begin{proposition}
Let $\{\Lambda_n\}_{n\in\mathbb N}$ be an $E$-$g$-frame for $\mathcal H$ with $E$-$g$-frame bounds $A$ and $B$ and the $E$-$g$-frame operator $S_{\Lambda}$. Then $\{\tilde{\Lambda}_n\}_{n\in\mathbb N}=\{\Lambda_nS_{\Lambda}^{-1}\}_{n\in\mathbb N}$ is an $E$-$g$-frame for $\mathcal H$.
\begin{proof}
For any $f\in\mathcal H$, we get
\begin{align*}
\sum_{n=1}^{\infty}\left\|\sum_{k=1}^{\infty}E_{n,k}\tilde{\Lambda}_kf\right\|^2&=
\sum_{n=1}^{\infty}\left\|\sum_{k=1}^{\infty}E_{n,k}\Lambda_kS_{\Lambda}^{-1}f\right\|^2\\
&=\sum_{n=1}^{\infty}\left\langle \sum_{k=1}^{\infty}E_{n,k}\Lambda_kS_{\Lambda}^{-1}f,\sum_{j=1}^{\infty}E_{n,j}\Lambda_jS_{\Lambda}^{-1}f \right\rangle\\
&=\left\langle\sum_{n=1}^{\infty}\sum_{k=1}^{\infty}\sum_{j=1}^{\infty}\overline{E}_{n,k}E_{n,j}\Lambda_j^*\Lambda_kS_{\Lambda}^{-1}f  ,S_{\Lambda}^{-1}f\right\rangle\\
&=\left\langle S_{\Lambda}S_{\Lambda}^{-1}f,S_{\Lambda}^{-1}f \right\rangle\\
&=\left\langle f,S_{\Lambda}^{-1}f \right\rangle\\
&\leq\frac{1}{A}\|f\|^2.
\end{align*}
On the other hand, we have
\begin{align*}
\|f\|^2&=\sum_{n=1}^{\infty}\sum_{k=1}^{\infty}\sum_{j=1}^{\infty}\left\langle\overline{E}_{n,k}E_{n,j}\Lambda_k^*\tilde{\Lambda}_jf,f\right\rangle\\
&=\sum_{n=1}^{\infty}\left\langle\sum_{j=1}^{\infty}E_{n,j}\tilde{\Lambda}_jf,\sum_{k=1}^{\infty}E_{n,k}\Lambda_kf\right\rangle\\
&\leq\left(\sum_{n=1}^{\infty}\left\|\sum_{j=1}^{\infty}E_{n,j}\tilde{\Lambda}_jf\right\|^2 \right)^{\frac{1}{2}}\left(\sum_{n=1}^{\infty}\left\|\sum_{k=1}^{\infty}E_{n,k}k\Lambda_kf\right\|^2 \right)^{\frac{1}{2}}\\
&\leq B^{\frac{1}{2}}\left(\sum_{n=1}^{\infty}\left\|\sum_{j=1}^{\infty}E_{n,j}\tilde{\Lambda}_jf\right\|^2 \right)^{\frac{1}{2}},
\end{align*}
this implies that
\begin{equation*}
\sum_{n=1}^{\infty}\left\|\sum_{j=1}^{\infty}E_{n,j}\tilde{\Lambda}_jf\right\|^2\geq\frac{1}{B}\|f\|^2.
\end{equation*}
Therefore, $\{\tilde{\Lambda}_n\}_{n\in\mathbb N}$ is an $E$-$g$-frame for $\mathcal H$ with $E$-$g$-frame bounds $\frac{1}{A}$ and $\frac{1}{B}$.
\end{proof}
\end{proposition}
\begin{remark}
Let $\tilde{S}_{\Lambda}$ be the $E$-$g$-frame operator associated with $\{\tilde{\Lambda}_n\}_{n\in\mathbb N}$. Then
\begin{align*}
S_{\Lambda}\tilde{S}_{\Lambda}&=\sum_{n=1}^{\infty}\sum_{j=1}^{\infty}E_{n,j}S_{\Lambda}\tilde{\Lambda}_k^*\tilde{\Lambda}_jf\\
&=\sum_{n=1}^{\infty}\sum_{j=1}^{\infty}E_{n,j}S_{\Lambda}S_{\Lambda}^{-1}\Lambda_k^*\Lambda_jS_{\Lambda}^{-1}f\\
&=\sum_{n=1}^{\infty}\sum_{j=1}^{\infty}E_{n,j}\Lambda_k^*\Lambda_jS_{\Lambda}^{-1}f\\
&=S_{\Lambda}S_{\Lambda}^{-1}f=f,
\end{align*}
for all $f\in\mathcal H$. Hence $\tilde{S}_{\Lambda}=S_{\Lambda}^{-1}$ and so $\tilde{\Lambda}_n\tilde{S}_{\Lambda}^{-1}=\Lambda_nS_{\Lambda}^{-1}S_{\Lambda}=\Lambda_n$. This asserts that $\{\Lambda_n\}_{n\in\mathbb N}$ and $\{\tilde{\Lambda}_n\}_{n\in\mathbb N}$ are dual $E$-$g$-frames with respect to each other.
\end{remark}
Recall that a sequence $\{a_n\}_{n=1}^{\infty}$ in $\mathbb C$ is said to be positively confined sequence if
\begin{equation*}
0<\inf_{n\in\mathbb N}|a_n|\leq \sup_{n\in\mathbb N}|a_n|<\infty.
\end{equation*}
Now, we give a result concerning the perturbation of $E$-$g$-frames.
\begin{theorem}\label{th2.1}
Let $\{\Lambda_n\}_{n\in\mathbb N}$ be an $E$-$g$-frame and $\Gamma_n\in\mathcal B(\mathcal H,\mathcal H_n)$, for all $n\in\mathbb N$. If there exist constants $0\leq \alpha ,\beta<\frac{1}{2}$, such that for every $f\in\mathcal H$,
\begin{align*}
\sum_{n=1}^{\infty}\left\|a_n\sum_{k=1}^{\infty}E_{n,k}\Lambda_kf-b_n\sum_{k=1}^{\infty}E_{n,k}\Gamma\right\|^2\leq\alpha\sum_{n=1}^{\infty}\left\|a_n\sum_{k=1}^{\infty}E_{n,k}\Lambda_kf\right\|^2+\beta\sum_{n=1}^{\infty}\left\|b_n\sum_{k=1}^{\infty}E_{n,k}\Gamma_kf\right\|^2.
\end{align*}
Then $\{\Gamma_n\}_{n\in\mathbb N}$ is a $E$–$g$-frame, where $\{a_n\}_{n=1}^{\infty}$ and $\{b_n\}_{n=1}^{\infty}$ are positively confined sequences.
\begin{proof}
Let $\{\Lambda_n\}_{n\in\mathbb N}$ be an $E$-$g$-frame with frame bounds $A$ and $B$. Then for all $f\in \mathcal H$ we get
\begin{align*}
\sum_{n=1}^{\infty}\left\|b_n\sum_{k=1}^{\infty}E_{n,k}\Gamma_kf\right\|^2&=\sum_{n=1}^{\infty}\left\|b_n\sum_{k=1}^{\infty}E_{n,k}\Gamma_kf\pm a_n\sum_{k=1}^{\infty}E_{n,k}\Lambda_kf\right\|^2\\
&\leq 2\sum_{n=1}^{\infty}\left\|b_n\sum_{k=1}^{\infty}E_{n,k}\Gamma_kf-a_n\sum_{k=1}^{\infty}E_{n,k}\Lambda_kf\right\|^2\\&+2\sum_{n=1}^{\infty}\left\|a_n\sum_{k=1}^{\infty}E_{n,k}\Lambda_kf\right\|^2\\
&\leq 2\left(\alpha\sum_{n=1}^{\infty}\left\|a_n\sum_{k=1}^{\infty}E_{n,k}\Lambda_kf\right\|^2+\beta\sum_{n=1}^{\infty}\left\|b_n\sum_{k=1}^{\infty}E_{n,k}\Gamma_kf\right\|^2\right)\\&\;+2\sum_{n=1}^{\infty}\left\|a_n\sum_{k=1}^{\infty}E_{n,k}\Lambda_kf\right\|^2.
\end{align*}
Therefore
\begin{align*}
(1-2\beta)\sum_{n=1}^{\infty}|b_n|^2\left\|\sum_{k=1}^{\infty}E_{n,k}\Gamma_kf\right\|^2\leq 2(1+\alpha)\sum_{n=1}^{\infty}|a_n|^2\left\|\sum_{k=1}^{\infty}E_{n,k}\Lambda_kf\right\|^2,
\end{align*}
this implies that
\begin{align*}
\sum_{n=1}^{\infty}\left\|\sum_{k=1}^{\infty}E_{n,k}\Gamma_kf\right\|^2&\leq\frac{2(1+\alpha)\left(\sup_{n\in\mathbb N}|a_n|\right)^2}{2(1-2\beta)\left(\inf_{n\in\mathbb N}|b_n|\right)^2}\sum_{n=1}^{\infty}\left\|\sum_{k=1}^{\infty}E_{n,k}\Lambda_kf\right\|^2\\
&\leq\frac{2(1+\alpha)\left(\sup_{n\in\mathbb N}|a_n|\right)^2}{2(1-2\beta)\left(\inf_{n\in\mathbb N}|b_n|\right)^2}B\left\|f\right\|^2.
\end{align*}
Similarly, we get
\begin{align*}
\sum_{n=1}^{\infty}\left\|a_n\sum_{k=1}^{\infty}E_{n,k}\Lambda_kf\right\|^2&\leq 2\alpha\sum_{n=1}^{\infty}\left\|a_n\sum_{k=1}^{\infty}E_{n,k}\Lambda_kf\right\|^2\\&+2\left(\beta\sum_{n=1}^{\infty}\left\|b_n\sum_{k=1}^{\infty}E_{n,k}\Gamma_kf\right\|^2+\sum_{n=1}^{\infty}\left\|b_n\sum_{k=1}^{\infty}E_{n,k}\Gamma_kf\right\|^2\right),
\end{align*}
consequently
\begin{align*}
\sum_{n=1}^{\infty}\left\|\sum_{k=1}^{\infty}E_{n,k}\Gamma_kf\right\|^2&\geq\frac{2(1-2\alpha)\left(\inf_{n\in\mathbb N}|a_n|\right)^2}{2(1+\beta)\left(\sup_{n\in\mathbb N}|b_n|\right)^2}\sum_{n=1}^{\infty}\left\|\sum_{k=1}^{\infty}E_{n,k}\Lambda_kf\right\|^2\\
&\geq\frac{2(1-2\alpha)\left(\inf_{n\in\mathbb N}|a_n|\right)^2}{2(1+\beta)\left(\sup_{n\in\mathbb N}|b_n|\right)^2}A\left\|f\right\|^2.
\end{align*}
Therefore $\{\Gamma_n\}_{n\in\mathbb N}$ is a $E$–$g$-frame.
\end{proof}
\end{theorem}
If we set $a_n=b_n=1$ for all $n$ and $\beta=0$ in Theorem \ref{th2.1}, then we have the following result.
\begin{corollary}\label{cor1}
Let $\{\Lambda_n\}_{n\in\mathbb N}$ be an $E$-$g$-frame and $\Gamma_n\in\mathcal B(\mathcal H,\mathcal H_n)$, for all $n\in\mathbb N$. 
If there exists $0\leq \alpha<\frac{1}{2}$, such that for every $f\in\mathcal H$,
\begin{align*}
\sum_{n=1}^{\infty}\left\|\sum_{k=1}^{\infty}E_{n,k}\Lambda_kf-\sum_{k=1}^{\infty}E_{n,k}\Gamma\right\|^2\leq\alpha\sum_{n=1}^{\infty}\left\|\sum_{k=1}^{\infty}E_{n,k}\Lambda_kf\right\|^2,
\end{align*}
Then $\{\Gamma_n\}_{n\in\mathbb N}$ is a $E$–$g$-frame.
\end{corollary}
Let $U\in \mathcal B(\mathcal H)$ and $M$ be a subspace of $\mathcal H$. Then, $M$ is said to be an
invariant subspace for $U$ if $Uf\in M$ whenever $f\in M$.\\
Since for every $f\in\mathcal H$ and $n\in\mathbb N$, we have
\begin{equation*}
\left\|\sum_{n=1}^{\infty}E_{n,k}\Lambda_kf-\left(\sum_{k=1}^{\infty}E_{n,k}\Lambda_k+U^mE_{n,k}\Lambda_k\right)f\right\|^2\leq\left\|U\right\|^{2m}\left\|E_{n,k}\Lambda_kf\right\|^2,
\end{equation*}
using Corollary \ref{cor1}, we get the following result.
\begin{corollary}
Let $U\in \mathcal B(\mathcal H)$ with $\|U\|<\frac{\sqrt{2}}{2}$, and $\mathcal H_n$ be invariant for $U$, and let $\{\Lambda_n\}_{n\in\mathbb N}\in\mathcal B(\mathcal H,\mathcal H_n)$, for all $n\in\mathbb N$. The followings are equivalent.
\begin{itemize}
\item[(i)]
$\{\Lambda_n\}_{n\in\mathbb N}$ is an $E$-$g$-frame.
\item[(ii)]
For each $m\in\mathbb N$, the sequence $\{\Lambda_n+U^m\Lambda_n\}_{n\in\mathbb N}$ is an $E$-$g$-frame.
\item[(iii)]
There exists $m\in\mathbb N$ such that $\{\Lambda_n+U^m\Lambda_n\}_{n\in\mathbb N}$ is an $E$-$g$-frame.
\end{itemize}
\end{corollary}
In the following theorem, we show that the product of an $E$-$g$-frame and some bounded linear operators is an $E$-$g$-frame.
\begin{theorem}
Let $U\in\mathcal B(\mathcal H)$ be a be a self-adjoint injective operator. Then $\{\Lambda_n\}_{n\in\mathbb N}$ is an $E$-$g$-frame if and only if  $\{\Lambda_nU\}_{n\in\mathbb N}$ is an $E$-$g$-frame. Also, for all $m\in\mathbb N$, the sequence $\{\Lambda_nU^m\}_{n\in\mathbb N}$ is an $E$-$g$-frame.
\begin{proof}
Let $U\in\mathcal B(\mathcal H)$ be a be injective, then there exists $U^{-1}\in\mathcal B(\mathcal H)$ such that $UU^{-1}=I_{\mathcal H}$. Hence for every $f\in \mathcal H$, we get
\begin{align*}
\|f\|&=\left\|(U^{-1})^*U^*f\right\|\\
&\leq\left\|(U^{-1})^*\right\|\left\|U^*f\right\|.
\end{align*}
This implies that
\begin{equation}\label{eq2.1}
\left\|(U^{-1})^*\right\|^{-1}\|f\|\leq\left\|U^*f\right\|.
\end{equation}
If $\{\Lambda_n\}_{n\in\mathbb N}$ is an $E$-$g$-frame, then  there exists $A>0$, such that 
for every $f\in \mathcal H$
\begin{align*}
\left\|\sum_{k=1}^{\infty}E_{n,k}\Lambda_kUf\right\|^2&\geq A\|Uf\|^2\\
&=A\left\|U^*f\right\|^2\\
&\geq A\left\|(U^{-1})^*\right\|^{-2}\left\|f\right\|^2.
\end{align*}
Hence $\{\Lambda_nU\}_{n\in\mathbb N}$ is an $E$-$g$-frame. 

Now, let $\{\Lambda_nU\}_{n\in\mathbb N}$ is an $E$-$g$-frame. Since $U$ is invertible and $U^{-1}$ is self-adjoint, hence $\{\Lambda_n\}_{n\in\mathbb N}=\{\Lambda_nUU^{-1}\}_{n\in\mathbb N}$ is an $E$-$g$-frame.
\end{proof}
\end{theorem}
Since the $E$-$g$-frame operator is self-adjoint and injective, we get immediately the following result.
\begin{corollary}
Let $U\in\mathcal B(\mathcal H)$ and $\{\Lambda_nU\}_{n\in\mathbb N}$ be an $E$-$g$-frame with the $E$-$g$-frame operator $S$. Then $\{\Lambda_nS^{-\frac{1}{2}}\}_n$ is an $E$-$g$-frame.
\end{corollary}
\begin{theorem}
Let $U\in\mathcal B(\mathcal H)$ and $\{\Lambda_nU\}_{n\in\mathbb N}$ be an $E$-$g$-frame. Then, the followings hold.
\begin{itemize}
\item[(i)]
$U$ is injective.
\item[(ii)]
If $U$ is self-adjoint and has closed range, then $U$ is invertible and $\{\Lambda_n\}_{n\in\mathbb N}$ is an $E$-$g$-frame.
\end{itemize}
\begin{proof}
\begin{itemize}
\item[(i)] 
It is clear by (\ref{def1}).
\item[(ii)]
Since $U$ is self-adjoint and has closed range, so
\begin{equation*}
R(U)=R(U^*)=N(U)^{\perp}=\mathcal H.
\end{equation*}
Therefore $U$ is surjective and by $(i)$ $U$ is invertible. 

Now, let $f\in\mathcal H$, so there exists $g\in\mathcal H$ such that $Ug=f$. In view of (\ref{eq2.1}), we get
\begin{align*}
\sum_{n=1}^{\infty}\left\|\sum_{k=1}^{\infty}E_{n,k}\Lambda_kf\right\|^2&=
\sum_{n=1}^{\infty}\left\|\sum_{k=1}^{\infty}E_{n,k}\Lambda_kUg\right\|^2\\
&\geq A\|g\|^2\\
&=A\left\|U^{-1}f\right\|^2\\
&=A\left\|(U^{-1})^*f\right\|^2\\
&\geq A\left\|U^*\right\|^{-2}\|f\|^2\\
&=A\left\|U\right\|^{-2}\|f\|^2.
\end{align*}
Therefore $\{\Lambda_n\}_{n\in\mathbb N}$ is an $E$-$g$-frame.
\end{itemize}
\end{proof}
\end{theorem}
In the next result, we give conditions under which the sum of an $E$-$g$-frame and an $E$-$g$-Bessel sequence is an $E$-$g$-Bessel sequence.
\begin{theorem}
Let $\{\Lambda_n\}_{n\in\mathbb N}$ be an 
$E$-$g$-frame with bounds $A_1$ and $B_1$ and let $\{\Gamma_n\}_{n\in\mathbb N}$ be an 
$E$-$g$-Bessel sequence with Bessel bound $B_2$. If for every $n\in\mathbb N$, $R(\Lambda_n)\perp  R(\Gamma_n)$, and 
$U_1,U_2\in B(\mathcal H)$, then 
$\{\Lambda_nU_1+\Gamma_nU_2\}_{n\in\mathbb N}$ is an $E$-$g$-Bessel sequence.
\begin{proof}
Since $R(\Lambda_n)\perp R(\Gamma_n)$, so 
$$\left\|\sum_{k=1}^{\infty}E_{n,k}\Lambda_kU_1f\right\|=\left\|\sum_{k=1}^{\infty}E_{n,k}\left(\Lambda_kU_1+\Gamma_kU_2\right)f\right\|,$$
for all $f\in\mathcal H$. Therefore
\begin{align*}
\sum_{n=1}^{\infty}\left\|\sum_{k=1}^{\infty}E_{n,k}\Lambda_kU_1f\right\|^2&\leq\sum_{n=1}^{\infty}\left\|\sum_{k=1}^{\infty}E_{n,k}\left(\Lambda_kU_1+\Gamma_kU_2\right)f\right\|^2\\
&\leq\sum_{n=1}^{\infty}\left\|\sum_{k=1}^{\infty}E_{n,k}\Lambda_kU_1f\right\|^2+\sum_{n=1}^{\infty}\left\|\sum_{k=1}^{\infty}E_{n,k}\Gamma_kU_2f\right\|^2\\
&\leq\left(B_1\|U_1\|^2+B_2\|U_2\|^2\right)\|f\|^2,
\end{align*}
for all $f\in\mathcal H$. Consequently, $\{\Lambda_nU_1+\Gamma_nU_2\}_{n\in\mathbb N}$ is an $E$-$g$-Bessel sequence.
\end{proof}
\end{theorem}
We get with the following result.
\begin{corollary}
Let $\{\Lambda_n\}_{n\in\mathbb N}$ be an $E$-$g$-frame with bounds $A_1$ and $B_1$ and let 
$\{\Gamma_n\}_{n\in\mathbb N}$ be an $E$-$g$-Bessel sequence with Bessel bound $B_2$. If for every $n\in\mathbb N$, $R(\Lambda_n)\perp R(\Gamma_n)$, then the followings hold.
\begin{itemize}
\item[(i)]
The sequences $\{\Lambda_n+\Gamma_n\}_{n\in\mathbb N}$ and $\{\Lambda_n-\Gamma_n\}_{n\in\mathbb N}$ are $E$-$g$-frames.
\item[(ii)]
If $\{a_n\}_{n=1}^{\infty}$ and $\{b_n\}_{n=1}^{\infty}$ are two positively confined sequences, then $\{a_n\Lambda_n+b_n\Gamma_n\}_{n\in\mathbb N}$ is an $E$-$g$-frame.
\end{itemize}
\end{corollary}

\end{document}